# Homotopy and Homology Vanishing Theorems and the Stability of Stochastic Flows [*]

K.D. Elworthy
Mathematics Institute
University of Warwick
Coventry CV4 7AL
UK

Steven Rosenberg
Mathematics Department
Boston University
Boston, MA 02215
USA

March 22, 1995

## §0. Introduction

A. In this paper we relate stability properties (*i.e.* moment exponents) of a stochastic dynamical system on a compact manifold $M$ to the homotopy and integral homology groups of $M$. In the special case of gradient Brownian systems associated to isometric immersions of $M$ in $\mathbb{R}^m$, these moment exponents can be estimated in terms of the second fundamental form of the immersion. This yields topological obstructions to isometric immersions generalizing results in [L-S],[H-W],[O], [Le] as well as new results on $p$-harmonic maps. At the same time, our work places these authors' results in the general framework of Weitzenböck formulas.

Recall that the Weitzenböck formula $\Delta^q = \nabla^*\nabla + R^q$ for the Laplacian on $q$-forms yields $H^q(M; \mathbb{R}) = 0$ provided the curvature term $R^q > 0$. This vanishing theorem has two drawbacks, namely that real cohomology contains limited topological information, and that the term $R^q$ is unmanagable for $q > 1$.

In contrast, it is shown in [L-S] that $H_q(M^n; \mathbb{Z}) = H_{n-q}(M^n; \mathbb{Z}) = 0$ provided

$$-\sum_{j=1}^{q} \sum_{\ell=q+t}^{n} (2|\beta(v_j, v_\ell)|^2 - \langle \beta(v_j, v_j),\ \beta(v_\ell, v_\ell)\rangle) + q(n-q) > 0, \qquad (0.1)$$

where $\beta$ is the second fundamental form for an isometric immersion of $M$ in $S^m$. This result is obtained by averaging the flow of a current in a given integral homology class over a finite dimensional family of gradient vector fields, and showing that on the average the mass of the current decreases provided (0.1) holds. Since every homology class has a minimal current, the homology must vanish.

---

[*]This work was partially supported by the NATO Collaborative Research Grants Programme and the SERC. The second author was partially supported by the NSF.

Roughly speaking, in this paper we average the flow of currents or spheres under the stochastic dynamical system. This is an infinite dimensional family parametrized by Wiener space, so it is not surprising that we obtain more information than for a finite dimensional average. In [E-R1], the authors produced vanishing theorems for real cohomology by dominating the heat flow on forms by the semigroup associated to a Schrödinger operator on functions, where the potential term was built from the curvature term $R^q$. In this paper we dominate the flow of a gradient SDS using a Schrödinger operator whose potential term involves both $R^q$ and the second fundamental form (Theorem 2B, Theorem 4B). Thus a positivity condition on the potential term gives homotopy and integral homology vanishing theorems.

For isometric immersions, we can check that our positivity condition for the potential term matches closely with (0.1). This is essentially why the infinite dimensional averaging extends their finite dimensional averaging results. Moreover, following [L-S] we need only consider the flow on primitive $p$-forms. For such forms the curvature term $R^q$ simplifies to a sum of sectional curvatures (generalizing the Ricci curvature).

Thus we to some extent overcome the two weaknesses of the Weitzenböck approach at the expense of working extrinsicly. The example of lens spaces $L$, which have $H_1(L; \mathbb{Z}) \neq 0$ but admit metrics of positive Ricci curvature, shows that some extrinsic condition is necessary, in this case a quantified version of our intuition that the lens space has nontrivial homology because it is more crumpled than its covering sphere.

B. To set the notation, on an $n$-dimensional compact smooth manifold $M$ consider a smooth vector field $Y$ and a smooth vector bundle map $X \colon \mathbb{R}^m \to TM$ of the trivial $m$-plane bundle over $M$ into the tangent bundle of $M$. The latter corresponds to $m$ vector fields $X_1, ..., X_m$ given by $X_i(x) = X(x, e_i)$, for $e_1, ..., e_m$ the standard basis for $\mathbb{R}^m$.

With $X$ vanishing identically the vector field $Y$ determines a dynamical system $dx_t = Y(x_t)dt$ whose behaviour is strongly limited by the topology of $M$, for example via Morse theory when $Y$ is a gradient. For general $X$ our data determine a second order semi-elliptic operator
$$\mathcal{A} \equiv \mathcal{L}_X^2 + Y \qquad (0.2)$$
acting on functions on $M$, where $\mathcal{L}_X^2$ refers to the sum of the repeated Lie derivatives $\mathcal{L}_X^2(f) = \sum_{i=1}^m \mathcal{L}_{X^i}\mathcal{L}_{X^i}(f) = \sum_{i=1}^m X^i(X^i f)$. This is elliptic if and only if each $X(x) \colon \mathbb{R}^m \to T_x M$, $x \in M$, is onto, in which case we will call $(X, Y)$ *non-degenerate*.

Our approach is to consider the stochastic dynamical system

$$dx_t = X(x_t) \circ dB_t + Y(x_t)dt, \tag{0.3}$$

where $\{B_t: t \geq 0\}$ is a Brownian motion on $\mathbb{R}^n$, whose solution is a Markov process with differential generator $\mathcal{A}$. If $\{\Omega, \mathcal{F}, \mathbb{P}\}$ is the probability space of $\{B_t: t \geq 0\}$ so that $B_t: \Omega \to \mathbb{R}^m$, equation (0.3) has a solution flow $F_t: \Omega \times M \to M$, $t \geq 0$, consisting of random smooth diffeomorphisms of $M$, continuous in $t$. As for dynamical systems, the most tractable case is of gradients: $Y = 0$, $X^i = \nabla f^i$, for $f = (f^1, \ldots, f^m): M \to \mathbb{R}^m$ an immersion, where each solution $F_t(x_0)$, $t \geq 0$ is a Brownian motion on $M$.

C. The organization and main results of the paper are as follows. In §§1,2 we define the strong moment exponent maps $\mu_M = \mu_M^1 : \mathbb{R} \to \mathbb{R}$, $\mu_M^q : \mathbb{R} \to \mathbb{R}$, $q = 1, 2, \ldots$ of (0.2). The moment exponent $\mu^q(p)$ is an $L^p$-measure of the long time behavior of the flow $F_t$ acting on $q$-forms. Roughly speaking, the stability condition $\mu^q(1) < 0$ implies an exponential decay to the flow and hence a shrinking of representatives of $q$-homology or homotopy classes. In particular, $\mu_M(1) < 0$ implies $\pi_1(M) = 0$, $\mu_M^2(1) < 0$ implies $\pi_2(M) = 0$, $\mu_M^q(1) < 0$ implies $H_q(M; \mathbb{Z}) = 0$, and consequently $\mu_M\left(\left[\frac{n+1}{2}\right]\right) < 0$ implies $M^n$ is a homotopy sphere (Corollaries 1B.1, 1C.1, 2B.1). For completeness we recall some of the basic results about integral currents in §2. In §3 conditions on general $(X, Y)$ are given which imply these various forms of moment stability. In §4 we consider the gradient Brownian systems given by isometric immersions of $M$ in $\mathbb{R}^m$. We show that $\mu_M^q(p) < 0$ if $\Delta - 2h_p^q > 0$, where $\Delta$ is the (positive definite) Laplacian on functions on $M$ and $h_p^q$ is a function on $M$ determined by the second fundamental form of the immersion. This is in essence a Feynman-Kac type argument. In particular we see that $\Delta - 2h_p^q > 0$ implies $H_q(M; \mathbb{Z}) = 0$. As mentioned above, this is similar to intrinsic results in [E-R1], [Li2], where it is shown that $\Delta + R_0^q > 0$ implies $H_q(M; \mathbb{R}) = 0$ and $\Delta + \text{Ric}_0 > 0$ implies $\pi_1(M)$ finite. Here $R_0^q : M \to \mathbb{R}$ takes $x \in M$ to the lowest eigenvalue of the endomorphism $R^q(x)$, and similarly for the Ricci curvature Ric. (Note that Li's result improves Myers' theorem that positive Ricci curvature implies $\pi_1(M)$ finite.)

§5 is devoted to simplifying the expression for $h_p^q$ to known geometric quantities, and involves no probability. The arguments are based on a rather simple expression for the Laplacian on $q$-forms $\phi$ in terms of an isometric embedding:

$$\Delta^q \phi = -\sum_{i=1}^m \mathcal{L}_{X^i} \mathcal{L}_{X^i} \phi, \tag{0.4}$$

for $X^i = \nabla f^i$ as above. (A new proof of (0.4) is given in the Appendix.) To state our main result in sharp form, set

$$\hat{R}_0^q(x) = \inf\{\langle R^q(x)V, V\rangle : V \text{ is a primitive vector in } T_x^*M\}.$$

Recall that $V$ is primitive if $V = \ell^1 \wedge \cdots \wedge \ell^q$ for orthonormal vectors $\ell^i \in T_x^*M$. Of course $\hat{R}_0^1 = R_0^1 = \mathrm{Ric}_0$.

**Theorem 5A:** *Let $M^n$ be isometrically immersed in $\mathbb{R}^m$ with second fundamental form $\alpha$ and mean curvature $H\colon M \to \mathbb{R}^m$. If $\|\alpha\|\colon M \to \mathbb{R}$ is the pointwise Hilbert-Schmidt norm of $\alpha$, then*

$$\Delta + \hat{R}_0^q - \frac{\|\alpha\|^2}{2} + \frac{n}{2}|H|^2 > 0$$

*implies $H_q(M;\mathbb{Z}) = H_{n-q}(M;\mathbb{Z}) = 0$. Moreover, if $\Delta + \mathrm{Ric}_0 - \frac{\|\alpha\|^2}{2} + \frac{n}{2}|H|^2 > 0$, then $\pi_1(M) = 0$, and if $\Delta + \hat{R}_0^2 - \frac{\|\alpha\|^2}{2} + \frac{n}{2}|H|^2 > 0$, then $\pi_2(M) = 0$.*

As a corollary, the homology groups vanish if $\hat{R}_0^q > \frac{\|\alpha\|^2}{2} - \frac{n}{2}|H|^2$ pointwise. If the manifold is isometrically immersed in the unit sphere, this corollary is equivalent to the main result of [L-S]. By work of [E-R], [R-Y], the probabilistic approach strengthens [L-S] to give the same vanishing provided the pointwise condition holds on "most" of $M$.

As mentioned, because $\hat{R}_0^q$ involves only primitive vectors, it is given by taking the infimum of a sum of sectional curvatures (Proposition 5B). As a result, we can conclude that $M$ is a homotopy sphere if the sectional curvatures are larger than $\|\alpha\|^2/(2n-2)$ (Corollary 5B).

D. We would like to acknowledge conversations with S. Deshmukh, J. Eells, Y. Le Jan, X.-M. Li, M. Micallef, and J. Rawnsley.

# 1 Flows of stochastic dynamical systems: homotopy obstructions to moment stability.

A. For the solution flow $\{F_t\colon t \geq 0\}$ of our stochastic differential equation (0.3) on $M$ there is the derivative flow $TF_t$ on the tangent bundle TM

$$T_{x_0}F_t(-,\omega)\colon T_{x_0}M \to T_{F_t(x_0,\omega)}M \qquad x_0 \in M,\, t \geq 0,\, \omega \in \Omega.$$

If $v_t(\omega) = T_{x_0}F_t(v_0,\omega)$, where $v_0 \in T_{x_0}M$, then $\{v_t\colon t \geq 0\}$ satisfies a certain stochastic differential equation on TM. However its behaviour is more concisely expressed by the covariant equation

$$Dv_t = \nabla X(v_t) \circ dB_t + \nabla Y(v_t)dt \tag{1.1}$$

along the paths $\{x_t: t \geq 0\}$ for $x_t(\omega) = F_t(x_0, \omega)$, using the Levi-Civita connection on $M$.

For $q = 0, 1, 2, \ldots, n = \dim M$ there is the induced map of $q$-vectors

$$\Lambda^q T F_t(-, \omega): \Lambda^q TM \to \Lambda^q TM \qquad \omega \in \Omega, \ t \geq 0$$

so that
$$\Lambda^q T F_t(V, \omega) = F_t(-, \omega)_*(V) \in \Lambda^q T_{x_t} M \qquad V \in \Lambda^q T_{x_0} M.$$

For $p \in \mathbb{R}$ define the *p-th (strong) moment exponent* $\mu_M(p)$ of (0.3) by

$$\mu_M(p) = \overline{\lim}_{t \to \infty} \frac{1}{t} \log \sup_{x \in M} \mathbb{E} \|T_x F_t\|^p \tag{1.2}$$

Here a Riemannian metric on $M$ is used to define the Hilbert-Schmidt norm $\|T_x F_t\|$ but the result is independent of the choice of metric, $M$ being compact. From now on we will choose some metric and use it and its Levi-Civita connection. For $q = 1, \ldots, n$ this can be extended to define

$$\mu_M^q(p) = \overline{\lim}_{t \to \infty} \frac{1}{t} \log \sup_{x \in M} \mathbb{E} |\Lambda^q T_x F_t|^p. \tag{1.3}$$

Clearly $\mu_M^1(p) = \mu_M(p)$ and in general, for $p \geq 1$.

$$\mu_M^q(p) \leq \mu_M(qp). \tag{1.4}$$

The corresponding exponents $\mu_x(p)$ and $\mu_x^q(p)$ obtained by removing the supremum over $x$ in (1.3) and (1.4) have been studied by various authors e.g. see [A] and [E2] together with its references, and in depth by Baxendale and Stroock in [B-S]. In the latter paper ergodicity assumptions are imposed via a hypoellipticity assumption on the differential generator of the processes $\{\frac{v_t}{|v_t|}: t \geq 0\}$ on the tangent sphere bundle of $M$. With such an assumption Baxendale and Stroock show that $\mu_x(p) = \mu_M(p)$ for all $x \in M$, $p \in \mathbb{R}$. Following [A], as in [E2] we have for each $q$,

$$p \mapsto \mu_M^q(p) \quad \text{is convex,}$$
$$p \mapsto \frac{1}{p} \mu_M^q(p) \quad \text{is increasing.}$$

Clearly $\mu_M^q(0) = 0$. We will say that our system is *strongly moment stable* if $\mu_M(1) < 0$ and *strongly p-th moment stable* if $\mu_M(p) < 0$. Topological and geometric consequences of such stability have been pointed out in [E3], see also [E-Y], [Li1]. Our emphasis on strong moment stability allows more detailed results for more general systems than those given in [E2], [E-Y].

B. We now relate moment stability of the general SDS (0.3) to the topology of $M$. For $p \in \mathbb{R}$ recall that the $p$-energy $E_p(\alpha)$ of a $C^1$ map $\alpha\colon N \to M$ of a Riemannian manifold $N$ into $M$ is given (up to constant multiples) by

$$E_p(\alpha) = \int_N \|T_y\alpha\|^p dy$$

(the usual energy, used for harmonic maps, being the 2-energy).

**Theorem 1B.** *If $M$ admits an equation (0.3) which is strongly $p$-th moment stable then every $C^1$ map $\alpha_0\colon N \to M$ of finite $p$-energy of a Riemannian manifold $N$ into $M$ (with any Riemannian metric) is homotopic to a map with arbitrarily small $p$-energy, as is every $C^1$ map $\beta_0\colon M \to N$.*

**Proof.** For such $\alpha_0$ consider the energy of the random map $\alpha_t(\omega)\colon N \to M$ given by

$$\alpha_t(\omega)(y) = F_t(\alpha_0(y), \omega)$$

and take its expectation:

$$\begin{aligned}
\mathbb{E}(E_p(\alpha_t)) &= \mathbb{E}\int_N \|TF_t \circ T_y\alpha_0\|^p dy = \int_N \mathbb{E}(\|TF_t \circ T_y\alpha_0\|^p) dy \\
&\leq \text{const.} \int_N (\mathbb{E}\|T_{\alpha_0(y)}F_t\|^p) \cdot \|T_y\alpha_0\|^p dy < \text{const.} \sup_{x \in M}(\mathbb{E}\|T_xF_t\|^p) E_p(\alpha_0) \\
&\to 0 \text{ as } t \to \infty.
\end{aligned}$$

The first assertion follows since each $\alpha_t(\omega)$ is homotopic to $\alpha_0$. For the second, set $\beta_t = \beta_0 \circ F_t$ and observe that $\mathbb{E}(E_p(\beta_t)) \leq (\|d\beta_0\|_\infty)^p \int_M \mathbb{E}\|T_xF_t\|^p dx$. □

**Corollary 1B.1** *(i) Strong moment stability can only occur when $M$ is simply connected.*
*(ii) Strong 2-moment stability implies $\pi_1(M) = 0$ and $\pi_2(M) = 0$.*

**Proof.** Part (i) comes from the existence of minimal length geodesic loops in any homotopy class and (ii) from [W]. □

A version of (i) when $M$ is non-compact is given in [Li1].

C. Corollary 1B.1 (ii) can be refined using $\mu_M^2(1)$. Indeed:

**Theorem 1C.** *If $\mu_M^q(1) < 0$ every $C^1$ map $\alpha_0\colon N \to M$ of a compact oriented $q$-dimensional manifold $N$ into $M$ is homotopic to one of arbitrarily small volume.*

**Proof.** Proceed as for Theorem 1B. With $\alpha_t = F_t \circ \alpha_0$ and a fixed Riemannian metric on $N$, the volume of $\alpha_t(N)$ satisfies

$$\operatorname{Vol} \alpha_t(N) \leq const. \int_N \|\Lambda^q T_y \alpha_t\| dy.$$

giving

$$\mathbb{E} \operatorname{Vol} \alpha_t(N) \leq \text{const.} \int_N \|\Lambda^q T_y \alpha_0\| dy \cdot \sup_{x \in M} \mathbb{E} \|\Lambda^q T_x F_t\| \to 0.$$

$\square$

**Corollary 1C.1.** *If $\mu_M^2(1) < 0$ then $\pi_2(M) = 0$.*

**Proof.** By Sacks and Uhlenbeck [S-U], if $\pi_2(M) \neq 0$, then there is a non-trivial immersed 2-sphere of minimal area in its homotopy class. $\square$

# 2 Integral currents and homology obstructions to moment stability

A. It is not clear how to give a direct extension of the results of §1 to obtain vanishing for $\pi_q(M)$ with $q > 2$. Instead we will follow Lawson and Simon [L-S] and use integral currents to obtain vanishing results for $H_q(M; \mathbb{Z})$. First we summarize some of the main facts about integral currents as presented in [L-S], see also [F], [Mo].

Let $\mathcal{H}^q$ denote $q$-dimensional Hausdorff measure on $M$. A Borel subset $S$ of $M$ is *q-rectifiable* if it is a countable union of $q$-dimensional $C^1$ submanifolds of $M$ union a set of $\mathcal{H}^q$-measure zero. For such an $S$ there is a tangent space $T_x S$ in $T_x M$, with dimension $q$, for $\mathcal{H}^q$ almost all $x$ in $S$. An *orientation* of $S$ is an $\mathcal{H}^q$-measurable $\xi \colon S \to \Lambda^q TM$ such that for $\mathcal{H}^q$ almost all $x$ in $S$ we have $\xi_x = a_1(x) \wedge \cdots \wedge a_q(x)$ for linearly independent $a_1(x), \ldots, a_q(x)$ in $T_x S$ and $|\xi_x| = 1$. Two oriented $q$-rectifiable sets will be identified if they are the same outside a set of $\mathcal{H}^q$-measure zero. Any such $(S, \xi)$ determines a linear map $S \colon \Omega^q M \to \mathbb{R}$ on the space $\Omega^q M$ of smooth $q$-forms on $M$, given by $S(\phi) = \int_S \phi(\xi_x) d\mathcal{H}^q(x)$. Giving $\Omega^q M$ the sup-norm topology determines a norm on $\mathcal{S} \in (\Omega^q M)^*$ which is written $M(\mathcal{S})$. It turns out that $M(\mathcal{S}) = \mathcal{H}^q(S)$, the Hausdorff measure of $S$.

Let $R_q(M)$ be the additive subgroup of $(\Omega^q M)^*$ generated by such $\mathcal{S}$. Another result is that each $\mathcal{S} \in R_p(M)$ has the form $\mathcal{S} = \sum_{n=1}^{\infty} n \mathcal{S}_n$ with $\{\mathcal{S}_n\}_{n=1}^{\infty}$ corresponding to a family $(S_n, \xi_n)$ of disjoint oriented $q$-rectifiable sets having $M(\mathcal{S}) = \sum_{n=1}^{\infty} \mathcal{H}^q(S_n) < \infty$.

For $\mathcal{S} \in \Omega^q(M)^*$ define $\partial \mathcal{S} \colon \Omega^{q-1}(M) \to \mathbb{R}$ by $\partial \mathcal{S}(\phi) = \mathcal{S}(d\phi)$. Then $\mathcal{S}$ is called an *integral $q$-current* if $\mathcal{S} \in R_q(M)$ and $\partial \mathcal{S} \in R_{q-1}(M)$. Let $I_q(M)$ denote the additive group of all integral $q$-currents and set $I_*(M) = \oplus_{q=0}^{\infty} I_q(M)$. This becomes a chain complex using the induced map $\partial \colon I_*(M) \to I_*(M)$. Let $H_q(I_*(M))$ be the $q$-th homology group of this complex.

Any smooth map $f \colon M \to N$ of compact manifolds induces a homomorphism $f_* \colon I_q(M) \to I_q(N)$ by $f_*(\mathcal{S})(\phi) = \mathcal{S}(f^*\phi)$, which commutes with $\partial$. This gives an induced homomorphism $f_* \colon H_*(I_*(M)) \to H_*(I_*(N))$, so we have a functor from the category of smooth compact Riemannian manifolds and smooth maps to abelian groups and homomorphism. A basic result of Federer and Fleming [F-F] is that *this is naturally equivalent to the corresponding singular homology functor with integer coefficients.* In particular for each $q \geq 0$ there is an isomorphism $H_q(I_*(M)) \cong H_q(M; \mathbb{Z})$.

A fundamental theorem of Federer and Fleming [F-F] (cf. [Mo, Ch. 5]) states that for each compact Riemannian $M$ and $q \geq 0$ *every $\alpha \in H_q(I_*(M))$ has a representative current $\mathcal{S} \in \alpha$ of least area*, i.e. $M(\mathcal{S}) \leq M(\mathcal{S}')$ for all $\mathcal{S}' \in \alpha$. (In fact, all we will need below is that a homology class containing integral currents of arbitrarily small positive mass must be the zero class, which follows from the proof of Theorem 5.7 in [Mo].)

From [F, p. 387] we see that if $\mathcal{S} = (S, \xi)$ is $q$-rectifiable and $h \colon M \to M$ is a diffeomorphism then $h_*(\mathcal{S}) = (h[S], \eta)$ where $\eta(y) = \dfrac{\Lambda^q(T_x h)(\xi(x))}{\|\Lambda^q(T_x h)(\xi(x))\|}$ for $y = h(x)$. In particular
$$M(h_*\mathcal{S}) = \int_M \|\Lambda^q(T_x h)(\xi(x))\|_{h(x)} d\mathcal{H}^q(x) \tag{2.1}$$
(using the convention that $\xi(x) = 0$ if $x \notin S$).

B.   We can now easily prove our main result:

**Theorem 2B.** *If the flow $\{F_t \colon t \geq 0\}$ of the stochastic differential equation (0.3) on $M$ has*
$$\underline{\lim}_{t \to \infty} \frac{1}{t} \operatorname{ess}_q \sup_{x \in M} \log \mathbb{E}\|\Lambda^q T_x F_t\| < 0,$$
*where $\operatorname{ess}_q$ sup refers to Hausdorff $q$-measure, then $H_q(M; \mathbb{Z}) = 0$.*

**Proof.** Fix a class $\alpha \in H_q(M; \mathbb{Z})$ and a representative integral $q$-current, $\mathcal{S}$. We have $\mathcal{S} = \sum_{n=1}^{\infty} n\mathcal{S}_n$, for $\mathcal{S}_n = (S_n, \xi_n)$ with $\sum_{n=1}^{\infty} n\mathcal{H}^q(S_n) < \infty$. By (2.1),

$$\begin{aligned} M((F_t)_*\mathcal{S}) &= \sum_{n=1}^{\infty} nM((F_t)_*\mathcal{S}_n) \\ &= \int_M \sum_{n=1}^{\infty} n\|\Lambda^q(T_xF_t)(\xi_n(x))\|_{F_t(x)} d\mathcal{H}^q(x) \end{aligned}$$

as functions of $\omega \in \Omega$. Therefore

$$\begin{aligned} \mathbb{E}M((F_t)_*\mathcal{S}) &= \int_M \sum_{n=1}^{\infty} n\mathbb{E}\|\Lambda^q T_xF_t(\xi_n(x))\|_{F_t(x)} d\mathcal{H}^q(x) \\ &\leq \int_M \sum_{n=1}^{\infty} n\chi_{S_n}(x)\mathbb{E}\|\Lambda^q T_xF_t\|_x d\mathcal{H}^q(x) \\ &\leq \operatorname{ess}_q \sup_{x \in M}(\mathbb{E}\|\Lambda^q T_xF_t\|_x) \sum_{n=1}^{\infty} n\mathcal{H}^q(S_n), \end{aligned}$$

where $\chi_{S_n}$ is the indicator function of $S_n$. Thus $\underline{\lim}_{t \to \infty} \mathbb{E}M((F_t)_*\mathcal{S}) = 0$ as $t \to \infty$. However each $F_t(-, \omega)_*(\mathcal{S})$ represents $\alpha$ by naturality of the equivalence of $H_*(I_*(M))$ with singular homology and the fact that with probability one $F_t(-, \omega)$ is homotopic to the identity. Therefore since $\alpha$ contains currents of arbitrarily small positive mass, we must have $\alpha = 0$. $\square$

**Corollary 2B1.** *If a compact manifold $M$ admits a strongly $q$-th moment stable stochastic dynamical system then $\pi_k(M) = 0$ for $k = 1, \ldots, q$. In particular, for $q \geq \dfrac{\dim M}{2}$, a $q$-th moment stable equation can only exist on homotopy spheres.*

**Proof.** Since $\mu_M(0) = 0$ and $\mu_M(p)$ is convex, $\mu_M(q) < 0$ implies $\mu_M(p) < 0$ for $p = 1, \ldots, q$. By (1.4), $\mu_M^p(1) < 0$ and so $\underline{\mu}_M^p(1) < 0$ for $p = 1, \ldots, q$. Thus by the preceding theorem, $H_p(M; \mathbb{Z}) = 0$ in this range. The conclusions now follow by a standard argument in algebraic topology. $\square$

C. Of course the exponential decay needed for moment stability is not necessary for these vanishing results, although it is implied by the conditions on the coefficients of (0.3) given below. It is also worth noting the following remarks on the samplewise behavior, valid for any diffeomorphism $h$ replacing $F_t$:

**Remarks 2C.** (i) *Let $\mathcal{S}$ be a minimal area integral $q$-current in $M$. Then $\mathcal{H}^q(x \in \operatorname{supp} \mathcal{S} : \|\Lambda^q T_xF_t\| \geq 1) > 0$.*

(ii) *For any compact $M$, $\mathcal{H}^q(x \in M\colon \|T_xF_t\| \geq 1) > 0$, for $0 \leq q \leq n$. In particular*
$$\mathbb{E}\sup_{x \in M}\|T_xF_t\| \geq 1.$$

**Proof.** For (i) observe that

$$0 \leq M((F_t)_*\mathcal{S}) - M(\mathcal{S}) \leq \int_M \sum_{n=1}^\infty n\chi_{S_n}(x)(\|\Lambda^q T_xF_t\| - 1)d\mathcal{H}^q(x).$$

For (ii), since $\int_M 1\, dx = \int_M \det T_xF_t\, dx$, we have $\mathcal{H}^n(x \in M\colon \|\Lambda^n T_xF_t\| \geq 1) > 0$, and $\|\Lambda^n T_xF_t\| \leq \|T_xF_t\|^n$.

## 3   Conditions on the coefficients

A.   The expectations in the moment exponent in the main Theorem 2B can be calculated fairly explicitly. First we do it in terms of the deterministic flows determined by $X$. Let $S_t^i\colon M \to M$ be the solution flows of the vector fields $X^i$, $i = 0, 1, \ldots, m$, taking $X^0 = Y$. It is immediate from Itô's formula, as given in [E1], [E2], and the integrability of $\|T_xF_t\|^p$, cf. [E1], that:

**Proposition 3A.** *If $V_0 \in \Lambda^q T_{x_0}M$ and $V_t(\omega) = \Lambda^q(TF_t(-,\omega))V_0$ for $\omega \in \Omega$, $t \geq 0$ then*

$$\mathbb{E}|V_t|^p = |V_0|^p + \int_0^t \mathbb{E}\frac{\partial}{\partial r}|\Lambda^q(TS_r^0)(V_s)|^p\bigg|_{r=0}ds + \frac{1}{2}\sum_{i=1}^m \int_0^t \mathbb{E}\frac{\partial^2}{\partial r^2}|\Lambda^q(TS_r^i)(V_s)|^p\bigg|_{r=0}ds.$$

Note that $\dfrac{\partial}{\partial r}|\Lambda^q(TS_r^i)(V_s)|^p = p\langle \Lambda^q(TS_r^i)(V_s), \dfrac{D}{\partial r}\Lambda^q(TS_r^i)(V_s)\rangle |\Lambda^q(TS_r^i)(V_s)|^{p-2}$
and so

$$\begin{aligned}\frac{\partial^2}{\partial r^2}|\Lambda^q(TS_r^i)(V_s)|^p\bigg|_{r=0} &= p(p-2)\langle V_s, (d\Lambda^q)(\nabla X^i)(V_s)\rangle^2|V_s|^{p-4} + p|(d\Lambda^q)(\nabla X^i)(V_s)|^2|V_s|^{p-2} \\ &\quad + p\langle V_s, \frac{D^2}{\partial r^2}\Lambda^q(TS_r^i)(V_s)|_{r=0}\rangle|V_s|^{p-2},\end{aligned}$$

where $(d\Lambda^q)(\nabla X^i)$ is defined to be linear and satisfy

$$(d\Lambda^q)(\nabla X^i)(v_1 \wedge \cdots \wedge v_q) = \nabla X^i(v_1) \wedge \cdots \wedge v_q + \cdots + v_1 \wedge \cdots \wedge \nabla X^i(v_q)$$

for $v_1, \ldots, v_q \in T_xM$, $x \in M$. Also,

$$\begin{aligned}
\left.\frac{D^2}{\partial r^2}\Lambda^q(TS_r^i)(V_s)\right|_{r=0} &= \left.\frac{D}{\partial r}d\Lambda^q(\nabla X^i)(\Lambda^q(TS_r^i)(V_s))\right|_{r=0} \\
&= d\Lambda^q(\nabla^2 X^i(X^i(x_s)))V_s + d\Lambda^q(\nabla X^i)(d\Lambda^q(\nabla X^i)(V_s)) \\
&= d\Lambda^q\left(\nabla(\nabla X^i(X^i(-))) + R(X^i(x_s),-)X^i(x_s)\right)(V_s) \\
&\quad + \delta^2\Lambda^q(\nabla X^i)(V_s),
\end{aligned} \quad (3.1)$$

where $\delta^2\Lambda^q(\nabla X^i) = d\Lambda^q(\nabla X^i) \circ d\Lambda^q(\nabla X^i) - d\Lambda^q(\nabla X^i \circ \nabla X^i)$, $R$ is the curvature tensor, and $x_s = F_s(x_0)\colon \Omega \to M$.

Although we will only use (3.1) below, for completeness we summarize the computations in the non-degenerate case. Here the metric on $M$ may be chosen to be the top order symbol of $\mathcal{A}$, in which case $\mathcal{A} = \frac{1}{2}\Delta + Y^X$ for $Y^X = Y + \frac{1}{2}\sum_{i=1}^m \nabla X^i(X^i)$. Then

$$\begin{aligned}
\mathbb{E}|V_t|^p &= |V_0|^p + p\int_0^t |V_s|^{p-2}\mathbb{E}\langle V_s, d\Lambda^q(\nabla Y^X)(V_s)\rangle ds \\
&\quad + \frac{1}{2}\sum_{i=1}^m p\int_0^t \mathbb{E}|d\Lambda^q(\nabla X^i)(V_s)|^2 |V_s|^{p-2}ds \\
&\quad + \frac{1}{2}p(p-2)\int_0^t \sum_{i=1}^m \mathbb{E}\langle V_s, d\Lambda^q(\nabla X^i)(V_s)\rangle^2 |V_s|^{p-4}ds \\
&\quad + \frac{1}{2}p\int_0^t \sum_{i=1}^m \mathbb{E}\langle V_s, \delta^2\Lambda^q(\nabla X^i)(V_s)\rangle ds \\
&\quad - \frac{1}{2}p\int_0^t \langle d\Lambda^q(\mathrm{Ric}^\#)V_s, V_s\rangle ds
\end{aligned}$$

where $\mathrm{Ric}^\#\colon TM \to TM$ comes from the Ricci curvature, $\mathrm{Ric}^\#(v) = -\sum_{j=1}^n R(e_j, v)e_j$ for $e_1, \ldots, e_n$ orthonormal in $T_xM$, $v \in T_xM$.

# 4 Gradient Brownian systems: immersed manifolds in $\mathbf{R}^m$

A. Let $f\colon M \to \mathbb{R}^m$ be an immersion. Give $M$ the induced metric and take $X^i = \nabla f^i$, $i = 1, \ldots, m$, and set $Y = 0$. Then $X(x)\colon \mathbb{R}^m \to T_xM$ can be considered as the orthogonal projection for each $x \in M$. The resulting stochastic dynamical systems has $\mathcal{A} = \frac{1}{2}\Delta$ so that the solution processes are Brownian motions on $M$. These are the *gradient Brownian systems* and have special properties [E2], [E-Y].

Let $\nu_x$ be the space of normal vectors to $M$ at $x$. We have the second fundamental form $\alpha_x\colon T_xM \times T_xM \to \nu_x$ and shape operator $A_x\colon T_xM \times \nu_x \to T_xM$ related by $\langle \alpha_x(v_1, v_2), w\rangle = \langle A_x(v_1, w), v_2\rangle$. If $Z(x)\colon \mathbb{R}^m \to \nu_x$ is the orthogonal projection then $\nabla X(v)e = A_x(v, Z(x)e)$, for $v \in T_xM$, $e \in \mathbb{R}^m$, and, as shown in the Appendix,

$$\sum_{i=1}^m \frac{D^2}{\partial r^2}\Lambda^q(TS_t^i)(V_0)\Big|_{r=0} = -(R_x^q)^*(V_0) \tag{4.1}$$

if $V_0 \in \Lambda^q T_xM$, where $R^q$ is the Weitzenböck curvature term $R_x^q\colon \Lambda^q T_x^*M \to \Lambda^q T_x^*M$ given by the Weitzenböck formula for the Laplacian on $q$-forms $\Delta^q$

$$R_x^q(\varphi_x) = (d+\delta)^2(\varphi)_x + \operatorname{trace}(\nabla^2\varphi)_x = \Delta^q(\varphi)_x - \nabla^*\nabla(\varphi)_x$$

for $\varphi$ a $C^\infty$ $q$-form [C-F-K-S]. Recall that $(R_x^1)^*(v) = \operatorname{Ric}^\#(v)$ for $v \in T_xM$.

For $V \in \Lambda^q T_xM$ set

$$\begin{aligned}
H_p^q(V, V) &= \sum_i \left[ |d\Lambda^q(A(-, Z^i(x)))V|^2 + (p-2)|V|^{-2}\langle V, d\Lambda^q(A(-, Z^i(x)))V\rangle^2 \right] \\
&\quad - \langle V, (R_x^q)^*V\rangle \\
&= \sum_i \left[ |d\Lambda^q(A(-, Z^i(x)))V|^2 + (p-2)|V|^{-2}\langle V, d\Lambda^q(A(-Z^i(x)))V\rangle^2 \right. \\
&\quad \left. + \langle V, \delta^2\Lambda^q(A(-, Z^i(x)))V\rangle \right] - \langle d\Lambda^q(\operatorname{Ric}^\#)V, V\rangle
\end{aligned} \tag{4.2}$$

using (A.2). Here $Z^i(x) = Z(x)(e_i)$. Then

**Lemma 4A.** *For a gradient Brownian system*

$$|V_t|^p = |V_0|^p + p\sum_i \int_0^t |V_s|^{p-2}\langle V_s, d\Lambda^q(A(-, Z^i(x_s)))V_s\rangle dB_s^i + (p/2)\int_0^t H_p^q(V_s, V_s)ds. \tag{4.3}$$

**Proof.** This follows from the Itô formula and (4.1). □

B. We now relate the moment exponents to the positivity of a Schrödinger operator on functions on $M$. It is crucial to note that the moment exponents only depend on the flow of primitive vectors in $\Lambda^*TM$, since such vectors form a basis of $\Lambda^*TM$ and are preserved under the flow. So let $\mathcal{P}_x^q$ denote the set of primitive vectors in $\Lambda^q T_xM$, and define the norm of the quadratic form $H_p^q$ on $\mathcal{P}_x$ by $h_p^q(x) = \sup\{\frac{p}{2}H_p^q(V,V)\colon V \in \mathcal{P}_x^q \text{ and } |V| = 1\}$. Let $\lambda(h_p^q)$ denote the infimum of the spectrum of $\Delta - 2h_p^q$, i.e. the smallest eigenvalue in the compact case under consideration. The following was pointed out to us by X.-M. Li, cf. [Li2]:

**Theorem 4B.** (X.-M. Li) *For a gradient Brownian system*

$$\mu_M^q(p) \leq -\frac{1}{2}\lambda(h_p^q).$$

**Proof.** Take an orthonormal base $\{E_x^j\}$ for $\Lambda^q T_x M$ consisting of elements in $\mathcal{P}_x^q$. We can use the $\ell_\infty$ norm to calculate $\mu_M^q(p)$, so that

$$\mu_M^q(p) = \overline{\lim}_{t\to\infty} \frac{1}{t} \sup_x \log \mathbb{E} \sup_j |\Lambda^q(TF_t)(E_x^j)|^p.$$

Set $V_t^j(x) = \Lambda^q(T_x F_t)(E_x^j)$ and $W_t^j(x) = V_t^j(x)/|V_t^j(x)|$. Note that $W_t^j(x)$ is primitive for all $t \geq 0$. Set

$$M_t^{j,q} = \sum_i \int_0^t \langle W_s^j(x), d\Lambda^q(A(-,Z^i(x_s)))W_s^j(x)\rangle dB_s^i$$

and

$$a_{p,t}^{j,q} = \frac{p}{2} \int_0^t H_p^q(W_s^j(x), W_s^j(x))ds.$$

Then, by (4.3),

$$|V_t^j(x)|^p = \mathcal{E}(pM_t^{j,q})e^{a_{p,t}^{j,q}}$$

where $\mathcal{E}(pM_t^{j,q})$ is the exponential martingale

$$\mathcal{E}(pM_t^{j,q}) = \exp[pM_t^{j,q} - \frac{1}{2}p^2 \langle M^{j,q}, M^{j,q}\rangle_t].$$

Thus

$$|V_t^j(x)|^p \leq \mathcal{E}(pM_t^{j,q})e^{\int_0^t h_p^q(x_s)ds}.$$

If we now apply the Girsanov-Maruyama theorem by setting

$$\tilde{B}_t = B_t - \sum_i p\left(\int_0^t \langle W_t^j(x), d\Lambda^q(A(-,Z^i(x_s)))W_t^j(x)\rangle ds\right) e_i$$

for $e_1, \ldots, e_m$ the standard base for $\mathbb{R}^m$, we see that $\{\tilde{B}_t: 0 \leq t \leq T\}$ is a Brownian motion on $\mathbb{R}^m$ with respect to the probability $\mathbb{Q} = \mathcal{E}(pM_t^{j,q})\mathbb{P}$, for any $T > 0$. However

$$\begin{aligned}
dx_t &= X(x_t) \circ d\tilde{B}_t + \sum_i pX^i(x_t)\langle W_t^j(x), d\Lambda^q(A(-,Z^i(x_t)))W_t^j(x)\rangle dt \\
&= X(x_t) \circ d\tilde{B}_t
\end{aligned}$$

so that $\{x_t: 0 \leq t \leq T\}$ has the same law under $\mathbb{Q}$ as under $\mathbb{P}$. Thus

$$\mathbb{E}\left|V_t^j(x)\right|^p \leq \mathbb{E}^{\mathbb{Q}} e^{\int_0^t h_p^q(x_s)ds} = \mathbb{E} e^{\int_0^t h_p^q(x_s)ds},$$

and the result follows from the standard consequence of the Feynman-Kac formula

$$\frac{1}{2}\lambda(h_p^q) = -\lim_{t \to \infty} \frac{1}{t} \log \mathbb{E} \ e^{\int_0^t h_p^q(x_s) ds}.$$

□

This allows us to relate homotopy and integral homology groups to the positivity of a Schrödinger operator on functions.

**Corollary 4B.1** $\Delta - 2h_1^q > 0$ *implies* $H_q(M; \mathbb{Z}) = 0$. *Also,* $\Delta - 2h_1^1 > 0$ *implies* $\pi_1(M) = 0$, *and* $\Delta - 2h_2^1 > 0$ *implies* $\pi_2(M) = 0$.

**Proof.** The first statement follows from Theorems 2B and 4B and the others from Corollary 1B.1 and Theorem 4B. □

In particular, we get vanishing theorems provided $h_p^q \leq 0$ provided the inequality is strict at at least one point. The conditions $h_1^q < 0$ are essentially those which appeared in [L-S] and $h_2^1 < 0$ in [H-W] (though using the second form of (4.2), not involving the Weitzenböck curvature). In §5 we examine these conditions in detail.

C. **Example.**

For the sphere $S^n(r)$ of radius $r$ in $\mathbb{R}^{n+1}$, $A(v, \frac{x}{r}) = -\frac{1}{r}v$ for $v \in T_x S^n(r)$ and, from (4.2), we see

$$h_p^q(x) = \frac{1}{2}pH_p^q(V,V) = \frac{1}{2r^2}pq(pq - n)$$

for $V \in \Lambda^q T_x M$ with $|V| = 1$. Thus $h_p^q < 0$ provided $pq < n$, showing that our homotopy and homology vanishing conditions are precise in this case.

# 5 Vanishing theorems for immersed manifolds

A. We now simplify formulas from §4A. By Corollary 4B.1, we have $H_q(M; \mathbb{Z}) = 0$ if $\Delta - 2h_1^q > 0$, where

$$2h_1^q = \sup \{-\sum_i \langle V, d\Lambda^1(A(-, Z^i))V \rangle^2 + \sum_i |d\Lambda^q(A(-, Z^i))V|^2 - \langle R^q V, V \rangle\} \quad (5.1)$$

for all $V = v_1 \wedge \cdots \wedge v_q$ with $\{v_i\}$ orthonormal. Here and from now on we shall ignore the distinction between $R^q$ and $(R^q)^*$. Now

$$-\sum_i \langle v_1 \wedge \cdots \wedge v_q, d\Lambda^q(A(-, Z^i)) v_1 \wedge \cdots \wedge v_q \rangle^2$$

$$= -\sum_i \langle v_1 \wedge \cdots \wedge v_q, \sum_j v_1 \wedge \cdots \wedge A(v_j, Z^i) \wedge \cdots \wedge v_q \rangle^2$$

$$= -\sum_i \left( \sum_{j=1}^q \langle v_j, A(v_j, Z^i) \rangle \right)^2 = -\sum_i \left( \sum_{j=1}^q \langle \alpha(v_j, v_j), Z^i \rangle \right)^2,$$

and

$$\sum_i |d\Lambda^q(A(-, Z^i))v|^2$$

$$= \sum_i \langle \sum_j v_1 \wedge \cdots \wedge A(v_j, Z^i) \wedge \cdots \wedge v_q, \sum_k v_1 \wedge \cdots \wedge A(v_k, Z^i) \wedge \cdots \wedge v_q \rangle$$

$$= \sum_i \{ \sum_{j=1}^q \langle A(v_j, Z^i), A(v_j, Z^i) \rangle - \sum_{\substack{j,k=1 \\ j \neq k}}^q \langle \alpha(v_j, v_k), Z^i \rangle^2$$

$$+ \sum_{\substack{j,k=1 \\ j \neq k}}^q \langle \alpha(v_j, v_j), Z^i \rangle \langle \alpha(v_k, v_k), Z^i \rangle \}.$$

Thus the right hand side of (5.1) equals the supremum of

$$-\sum_i \left( \sum_{j=1}^q \langle \alpha(v_j, v_j), Z^i \rangle \right)^2 + \sum_i \left( \sum_{j=1}^q \sum_{\ell=1}^n \langle \alpha(v_j, v_\ell), Z^i \rangle^2 \right)$$

$$- \sum_{\substack{j,k=1 \\ j \neq k}}^q |\alpha(v_j, v_k)|^2 + \sum_i \sum_{\substack{j,k=1 \\ j \neq k}}^q \langle \alpha(v_j, v_j), Z^i \rangle \langle \alpha(v_k, v_k), Z^i \rangle - \langle R^q V, V \rangle$$

$$= -\sum_{j=1}^q |\alpha(v_j, v_j)|^2 + \sum_{j=1}^q \sum_{\ell=1}^n |\alpha(v_j, v_\ell)|^2$$

$$- \sum_{\substack{j,k=1 \\ j \neq k}}^q |\alpha(v_j, v_k)|^2 - \langle R^q V, V \rangle$$

$$= -\sum_{j=1}^q \sum_{\ell=q+1}^n |\alpha(v_j, v_\ell)|^2 - \langle R^q V, V \rangle \tag{5.2}$$

As usual, we write $R^q > C$ at $x \in M$ if $\langle R^q V, V \rangle > C|V|^2$ for all $V \in \Lambda^q T_x M$, $V \neq 0$. Recall that $R_0^q = R_0^q(x)$ is the smallest such $C$; equivalently, $R_0^q(x)$ is the lowest eigenvalue of $R^q$ on $\Lambda^p T_x M$. Similarly we set $\hat{R}_0^q(x) = \inf\{\langle R^q V, V \rangle \colon V \in \mathcal{P}_x^q, |V| = 1\}$, so that $\hat{R}_0^q(x) \geq R_0^q(x)$. By (5.2), we see that $\Delta - 2h_1^q > 0$ if

$$\Delta + \hat{R}_0^q - \sup \sum_{j=1}^q \sum_{\ell=q+1}^n |\alpha(v_j, v_\ell)|^2 > 0, \tag{5.3}$$

where the supremum is taken over all orthonormal bases $\{v_1, \ldots, v_n\}$ in $T_x M$.

At this point we can give our main vanishing theorem relating the Weitzenböck term on primitive forms to the second fundamental form. To state the result, let $\mathcal{N} = \mathcal{N}(n, K, D, V)$ be the class of Riemannian manifolds

$$\mathcal{N}(n, K, D, V) = \{(M^n, g)\colon \operatorname{Ric} \geq K,\ \operatorname{diam}\ \leq D,\ \operatorname{vol} \geq V\}.$$

Also, for $f\colon M \to \mathbb{R}$, set $f_-(x) = \min\{f(x), 0\}$.

**Theorem 5A.** *Let $M^n$ be isometrically immersed in $\mathbb{R}^m$ with second fundamental form $\alpha$, and mean curvature $H\colon M \to \mathbb{R}^m$. If*

$$\Delta + \hat{R}_0^q - \frac{\|\alpha\|^2}{2} + \frac{n}{2}|H|^2 > 0,$$

*then $H_q(M;\mathbb{Z}) = H_{n-q}(M;\mathbb{Z}) = 0$. Moreover, if*

$$\Delta + \operatorname{Ric}_0 - \frac{\|\alpha\|^2}{2} + \frac{n}{2}|H|^2 > 0,$$

*then $\pi_1(M) = 0$, and if*

$$\Delta + \hat{R}_0^2 - \frac{\|\alpha\|^2}{2} + \frac{n}{2}|H|^2 > 0,$$

*then $\pi_2(M) = 0$. In particular, there exists a constant $C = C(\mathcal{N}) > 0$ such that if (i) $M \in \mathcal{N}$ and (ii) there exists $w > 0$ such that*

$$\|(\hat{R}_0^q - \frac{\|\alpha\|^2}{2} + \frac{n}{2}|H|^2 - w)_-\|_{n/2} \leq \min\{(2\operatorname{vol}(M))^{2/n}, w \cdot C^{-1}(\mathcal{N})\},$$

*then $H_q(M;\mathbb{Z}) = H_{n-q}(M;\mathbb{Z}) = 0$.*

**Proof.** We have

$$\begin{aligned}
\frac{\|\alpha\|^2}{2} &= \frac{1}{2} \sum_{j,\ell=1}^{n} |\alpha(v_j, v_\ell)|^2 \\
&\geq \sum_{j=1}^{q} \sum_{\ell=q+1}^{n} |\alpha(v_j, v_\ell)|^2 + \frac{1}{2} \sum_{j=1}^{n} |\alpha(v_j, v_j)|^2 \\
&\geq \sum_{j=1}^{q} \sum_{\ell=q+1}^{n} |\alpha(v_j, v_\ell)|^2 + \frac{1}{2n} |\sum_{j=1}^{n} \alpha(v_j, v_j)|^2 \\
&= \sum_{j=1}^{2} \sum_{\ell=q+1}^{n} |\alpha(v_j, v_\ell)|^2 + \frac{n}{2}|H|^2.
\end{aligned}$$

Thus by (5.1) and (5.3), $\Delta + \hat{R}_0^q - \frac{\|\alpha\|^2}{2} + \frac{n}{2}|H|^2 > 0$ implies $H_q(M;\mathbb{Z}) = 0$. Since $*R^q = R^{n-q}*$, where $*$ is the Hodge star operator, $R_0^q = R_0^{n-q}$. Thus under the hypothesis we also obtain $H_{n-q}(M;\mathbb{Z}) = 0$. The statements about $\pi_1(M)$ and $\pi_2(M)$ follow similarly. The last statement is a consequence of [R-Y, Thm. 2.2]. (The isoperimetric constant in Thm. 2.2 is bounded above by a constant $C = C(\mathcal{N})$ by an elementary argument using the Bishop-Gromov comparison theorem.)

**Corollary 5A.** *If (i) $\hat{R}^q \geq \frac{\|\alpha\|^2}{2} - \frac{n}{2}|H|^2$ pointwise and $\hat{R}_q > \frac{\|\alpha\|^2}{2} - \frac{n}{2}|H|^2$ at some point in $M$, then $H_q(M;\mathbb{Z}) = H_{n-q}(M;\mathbb{Z}) = 0$.*

*(ii) If $\mathrm{Ric} \geq \frac{\|\alpha\|^2}{2} - \frac{n}{2}|H|^2$ pointwise and $\mathrm{Ric} > \frac{\|\alpha\|^2}{2} - \frac{n}{2}|H|^2$ at some point in $M$, then $\pi_1(M) = 0$.*

*(iii) If $\hat{R}^2 \geq \frac{\|\alpha\|^2}{2} - \frac{n}{2}|H|^2$ pointwise and $\hat{R}^2 > \frac{\|\alpha\|^2}{2} - \frac{n}{2}|H|^2$ at some point in $M$, then $\pi_2(M) = 0$.*

**Remarks.**

(i) As in [L-S] we may conclude $H_q(M;G) = H_{n-q}(M;G)$ for any finitely generated abelian group $G$.

(ii) The last statement of Theorem 5A shows that $H_q(M;\mathbb{Z}) = H_{n-q}(M;\mathbb{Z}) = 0$ if $\hat{R}^q > \frac{\|\alpha\|^2}{2} - \frac{n}{2}|H|^2$ on "most" of $M$. For example, if this inequality holds except on a set $V_1$, where the inequality fails by at most a constant $F$, then there exists a constant $c = c(F, \mathcal{N}) > 0$ such that $\mathrm{vol}(V_1) < c$ implies the vanishing of the homology (cf. [E-R1]). Thus we need not assume the pointwise conditions of the Corollary.

(iii) Under the hypotheses of Theorem 5A, we in fact have $H_q(M';\mathbb{Z}) = H_{n-q}(M';\mathbb{Z}) = 0$ for all finite covers $\pi: M' \to M$. For if $i: M \to \mathbb{R}^N$ is an isometric immersion, so is $i \circ \pi: M' \to \mathbb{R}^N$ with respect to the pullback metric on $M'$, and of course $R_0^q$ and $\alpha$ are the same for $M$ and $M'$.

In particular, a manifold with $\hat{R}^q > \frac{\|\alpha\|^2}{2} - \frac{n}{2}|H|^2$ for all $q$ is a homotopy sphere; this is in contrast to the usual Bochner result that a manifold with $R^q > 0$ for all $q$ is a real homology sphere. Note also that the usual intrinsic Bochner type vanishing results ($R^q > 0 \Rightarrow H^q(M;\mathbb{R}) = H^{n-q}(M;\mathbb{R}) = 0$) apply equally well to a manifold and a finite quotient. This is not the case in Theorem 5A. For example, Theorem 5A gives no information about spheres, but shows that for any lens space, $\mathrm{Ric} < \frac{\|\alpha\|^2}{2} - \frac{n}{2}|H|^2$ at some point of any isometric immersion or there is equality everywhere.

B. In Theorem 5A, we only used $R^p > C$ on primitive vectors $V = v_1 \wedge \cdots \wedge v_q$ with $|v_i| = 1$. We now obtain a simple expression for $\langle R^p V, V \rangle$ for such vectors. By equation (A.2) in the Appendix and then by Gauss' equation, after extending

$v_1, \ldots, v_q$ to an orthonormal base $v_1, \ldots, v_n$:

$$
\begin{aligned}
\langle -R^q V, V \rangle &= \langle V, \sum_i \delta^2 \Lambda^q(A(-, Z^i))V \rangle - \langle V, d\Lambda^q(\text{Ric}^\#)V \rangle \\
&= \sum_i \langle v_1 \wedge \cdots \wedge v_q, \sum_{\substack{j,k=1 \\ j \neq k}}^q v_1 \wedge \cdots A(v_j, Z^i) \wedge \cdots A(v_k, Z^i) \wedge \cdots \wedge v_q \rangle \\
&\quad - \sum_{j=1}^q \sum_{\ell=1}^n \langle \alpha(v_j, v_j), \alpha(v_\ell, v_\ell) \rangle + \sum_{j=1}^q \sum_{\ell=1}^n |\alpha(v_j, v_\ell)|^2 \\
&= \sum_i \sum_{\substack{j,k=1 \\ j \neq k}}^q \langle A(v_j, Z^i), v_j \rangle \langle A(v_k, Z^i), v_k \rangle \\
&\quad - \sum_i \sum_{\substack{j,k=1 \\ j \neq k}}^q \langle A(v_j, Z^i), v_k \rangle \langle A(v_k, Z^i), v_j \rangle \\
&\quad - \sum_{j=1}^q \sum_{\ell=1}^n \langle \alpha(v_j, v_j), \alpha(v_\ell, v_\ell) \rangle + \sum_{j=1}^q \sum_{\ell=1}^n |\alpha(v_j, v_\ell)|^2 \\
&= \sum_{\substack{j,k=1 \\ j \neq k}}^q \langle \alpha(v_j, v_j), \alpha(v_k, v_k) \rangle - \sum_{\substack{j,k=1 \\ j \neq k}}^q |\alpha(v_j, v_k)|^2 \\
&\quad - \sum_{j=1}^q \sum_{\ell=1}^n \langle \alpha(v_j, v_j), \alpha(v_\ell, v_\ell) \rangle + \sum_{j=1}^q \sum_{\ell=1}^n |\alpha(v_j, v_\ell)|^2 \\
&= \sum_{j=1}^q \sum_{k=1}^q \langle \alpha(v_j, v_j), \alpha(v_k, v_k) \rangle - \sum_{j=1}^q \sum_{k=1}^q |\alpha(v_j, v_k)|^2 \\
&\quad - \sum_{j=1}^q \sum_{\ell=1}^n \langle \alpha(v_j, v_j), \alpha(v_\ell, v_\ell) \rangle + \sum_{j=1}^q \sum_{\ell=1}^n |\alpha(v_j, v_\ell)|^2 \\
&= \sum_{j=1}^q \sum_{\ell=q+1}^n \left( -\langle \alpha(v_j, v_j), \alpha(v_\ell, v_\ell) \rangle + |\alpha(v_j, v_\ell)|^2 \right). \quad (5.4)
\end{aligned}
$$

From (5.4) and (5.2) we see that our basic criterion for strong $q$-moment stability, with consequent vanishing results, holds whenever

$$\sum_{j=1}^q \sum_{\ell=q+1}^n \{ \langle \alpha(v_j, v_j), \alpha(v_\ell, v_\ell) \rangle - 2|\alpha(v_j, v_\ell)|^2 \} \geq 0$$

for all orthonormal $v_1, \ldots, v_n$ in $T_x M$ for all points $x$ of $M$, with strict inequality at some point of $M$. This is essentially the criterion used in [L-S] and [H-W], see below, though they insist on strict inequality everywhere. Examples in [L-S] show that the vanishing does not follow without strict inequality somewhere, so these results are sharp.

By the Gauss equation, $-|\alpha(v_j, v_\ell)|^2 + \langle \alpha(v_j, v_j), \alpha(v_\ell, v_\ell) \rangle$ equals $K(v_j, v_\ell)$, the sectional curvature of the $v_j, v_\ell$ plane. This with (5.4) gives an intrinsic expression for $R^q$ on primitive vectors.

**Proposition 5B.** *Let $V = v_1 \wedge \cdots \wedge v_q \in \Lambda^q T_x M$ with $|v_i| = 1$. Then*

$$\langle R^q V, V \rangle = \sum_{j=1}^{q} \sum_{\ell=q+1}^{n} K(v_j, v_\ell),$$

*where $v_{q+1}, \ldots, v_n$ are chosen so that $v_1, \ldots, v_n$ is an orthonormal base for $T_x M$.*

An intrinsic proof can be obtained easily from the fermionic calculus formula

$$R_q = R_{ijk\ell}(v_i)^* v_j (v_k)^* v_\ell$$

for an orthonormal base $v_1, \ldots, v_n$ of tangent vectors with $v_j$ also denoting the operation of interior multiplication by $v_j$ and $(v_j)^*$ its adjoint (i.e. exterior multiplication by $v_j$), [C-F-K-S].

C. The curvature operator $\mathcal{R}: \Lambda^2 T_x^* M \to \Lambda^2 T_x^* M$ is defined by $\langle \mathcal{R}(X \wedge Y), Z \wedge W \rangle = \langle R(X, Y)Z, W \rangle$, when $R$ is the curvature tensor and we do not distinguish between a tangent vector and its dual. In particular, $\langle \mathcal{R}(X \wedge Y), X \wedge Y \rangle = K(X, Y)$ if $|X| = |Y| = 1$. By Proposition 5B, if $\mathcal{R} > C$ on primitive vectors, then $R^q > q(n-q)C$ on primitive vectors; the same statement with "primitive" omitted is in [G-M]. As before, $\hat{\mathcal{R}}_0(x)$ denotes the largest $C$ such that $\mathcal{R} > C$ on primitive vectors, so $\hat{\mathcal{R}}_0(x) = \inf\{K(X, Y) \colon X, Y \in T_x M\}$.

**Corollary 5C** (i). *If $\Delta + (n-1) \hat{\mathcal{R}}_0 - \frac{\|\alpha\|^2}{2} + \frac{n}{2}|H|^2 > 0$, then $M$ is a homotopy sphere. In particular if the sectional curvatures $K = K(X, Y)$ satisfy $K \geq \frac{\|\alpha\|^2}{2(n-1)} - \frac{n}{2(n-1)}|H|^2$ pointwise and $K > \frac{\|\alpha\|^2}{2(n-1)} - \frac{n}{2(n-1)}|H|^2$ at some point, then $M$ is a homotopy sphere.*

*(ii) If the sectional curvatures satisfy $K(u, v) \geq |\alpha(u, v)|^2$ whenever $u, v$ are orthogonal, for all points of $M$, with strict inequality at at least one point of $M$, then $M$ is a homotopy sphere.*

**Proof.** For (i) we just note

$$\Delta + \hat{R}_0^q - \frac{\|\alpha\|^2}{2} + \frac{n}{2}|H|^2 > \Delta + q(n-q)\hat{\mathcal{R}}_0 - \frac{\|\alpha\|^2}{2} + \frac{n}{2}|H|^2$$

$$\geq \Delta + (n-1)\hat{\mathcal{R}}_0 - \frac{\|\alpha\|^2}{2} + \frac{n}{2}|H|^2 > 0.$$

Part (ii) follows from equation (5.2) and the proposition. □

Leung [Le] proved (ii) under the assumption of strict inequality. The proof of (i) also shows that $K \geq \dfrac{\|\alpha\|^2}{2q(n-q)} - \dfrac{n}{2q(n-q)}|H|^2$ everywhere with strict inequality somewhere implies $H_q(M;\mathbb{Z}) = H_{n-q}(M;\mathbb{Z}) = 0$.

Note that the Gauss equation gives

$$\begin{aligned}\|\alpha\|^2 &= \sum_{j,\ell} |\alpha(v_j,v_\ell)|^2 = -\sum_{j,\ell} K(v_j,v_\ell) + |\sum_j \alpha(v_j,v_j)|^2 \\ &\leq -2\sum_{j<\ell} K(v_j,v_\ell) + n\sum_j |\alpha(v_j,v_j)|^2 \qquad (5.5) \\ &\leq -2\sum_{j<\ell} K(v_j,v_\ell) + n\|\alpha\|^2,\end{aligned}$$

and so the minimal sectional curvature always satisfies $K_{\min} \leq \|\alpha\|^2/n$. Thus Corollary 5C can be interpreted as showing that a pinching condition on $K_{\min}$ implies $M$ is a homotopy sphere.

Micallef and Moore [M-M] have shown that a simply connected manifold with $\mathcal{R}$ positive on complex isotopic two-planes is a homotopy sphere. In contrast, Corollary 5C replaces this pointwise intrinsic hypothesis with a global extrinsic one. Note also from (5.5) that $K_{\max} \leq (n-1)\dfrac{\|\alpha\|^2}{2}$, if $K > 0$ on $M$, so the pointwise hypothesis in Corollary 5B (i) only implies a pinching constant of $\dfrac{1}{(n-1)^2}$ for the sectional curvature. Thus the conclusion that $M$ is a homotopy sphere does not follow from standard sphere theorems, even if $M$ is simply connected and of positive sectional curvature.

As usual, the hypotheses on the sectional curvature in the Corollary can be weakened to holding on "most" of $M$. Finally, note the first line of (5.5) gives $\frac{1}{2}\|\alpha\|^2 - \frac{n}{2}|H|^2 = -\frac{1}{2}k + \frac{1}{2}n(n-1)|H|^2$ where $k$ is the scalar curvature.

D. We now show that Theorem 5A is a strengthening of results in [L-S]. Assume that $M$ is isometrically immersed in $S^N$, the standard unit sphere. Let $\beta, \alpha, \alpha'$ denote the second fundamental forms of $M$ in $S^N$, $M$ in $\mathbb{R}^{N+1}$ under the inclusion $S^N \hookrightarrow \mathbb{R}^{N+1}$, and $S^N$ in $\mathbb{R}^{N+1}$, respectively. From $\alpha'(v_1,v_2) = -\langle v_1,v_2\rangle \nu$, where $\nu$ is the unit outward normal to $S^N$, we find

$$\beta(v_1,v_2) - \langle v_1,v_2\rangle \nu = \alpha(v_1,v_2) \qquad (5.6)$$

In particular, $|\beta(v_1,v_2)|^2 = |\alpha(v_1,v_2)|^2$ if $v_1 \perp v_2$, while $|\beta(v_1,v_1)|^2 + 1 = |\beta(v_1,v_1)|^2$ if $|v_1| = 1$. Combining (5.3), (5.4), (5.6) gives the vanishing of $H_q(M;\mathbb{Z})$ and $H_{n-q}(M;\mathbb{Z})$ provided

$$\sum_{j=1}^{q}\sum_{\ell=q+1}^{n} \left(-|\beta(v_j,v_\ell)|^2 + \langle \beta(v_j,v_j)-\nu,\, \beta(v_\ell,v_\ell)-\nu\rangle\right) > \sum_{j=1}^{q}\sum_{\ell=q+1}^{n} |\beta(v_j,v_\ell)|^2$$

This yields

**Theorem 5D.** (Lawson-Simons). *Let $M_n$ be isometrically immersed in $S^N$ with second fundamental form $\beta$. If*

$$\sum_{j=1}^{q} \sum_{\ell=q+t}^{n} (2|\beta(v_j, v_\ell)|^2 - \langle \beta(v_j, v_j), \beta(v_\ell, v_\ell) \rangle) < q(n-q)$$

*then $H_q(M; \mathbb{Z}) = H_{n-q}(M : \mathbb{Z}) = 0$.*

Of course, the hypothesis on $\beta$ need only hold on "most" of $M$, and the inequality $<$ may be replaced by $\leq$, with strict inequality somewhere. Consequences of this theorem in [L-S] admit similar improvements. Let $\hat{N}(M)$ be the unit normal bundle of $M$ in $S^N$ with the metric induced from $TS^N$. The volume, diameter and sectional curvatures of $\hat{N}(M)$ are bounded in terms of the volume and diameter of $M$ and $\beta$. In particular, the isoperimetric constant for $\hat{N}(M)$ is easily seen to be bounded above by a constant $C = C_1(\mathcal{N}, \overline{\beta})$, if $M \in \mathcal{N}$, where $\overline{\beta}$ is the supremum of the pointwise norms $\|\beta\|_x^2$ (cf. the proof of Theorem 5A.)

Let $\kappa = \kappa(x, \nu)$ be the minimum of the squares of the principal curvatures of the second fundamental form $\beta^*$ of the polar map $\hat{N}(M) \to S^N$ at the point $(x, \nu) \in \hat{N}(M)$. Then $\|\beta^*\|_{(x,\nu)}^2 \leq n/\kappa$. Since

$$\sum_{j=1}^{q} \sum_{\ell=q+1}^{n} (2|\beta^*(v_j, v_\ell)|^2 - \langle \beta^*(v_j, v_j), \beta^*(v_\ell, v_\ell) \rangle) \leq \max\{1, \frac{\sqrt{q(N-1-q)}}{2}\} \|\beta^*\|_{(x,\nu)}^2$$

[L-S, p. 441], the argument giving Theorem 5D shows that $H_q(\hat{N}(M); \mathbb{Z}_2) = H_{N-1-q}(\hat{N}(M); \mathbb{Z}_2) = 0$ provided

$$\Delta + \kappa - \max\left\{\frac{n}{q(N-1-q)}, \frac{n}{2\sqrt{q(N-1-q)}}\right\} < 0 \qquad q \neq 0, N-1 \qquad (5.7)$$

on $\hat{N}(M)$.

**Corollary 5D:** (i) *Let $M^n$ be isometrically embedded in $S^{n+k}$, for $k > 1$, with second fundamental form $\beta$. There exist constants $C_1 = C_1(\mathcal{N}, A)$, $C_2 = C_2(\mathcal{N}, A)$ such that if $M \in \mathcal{N} = \mathcal{N}(n, K, D, V)$ and $\overline{\beta} < A$, then for all $w > 0$*

$$\|(\kappa - \max\left\{\frac{1}{k-1}, \frac{1}{2}\sqrt{\frac{n}{k-1}}\right\} - w)_-\|_{\frac{N-1}{2}} \geq \min\{C_1(\mathcal{N}, A), w \cdot C_2(\mathcal{N}, A)\}.$$

(ii) *Either* $\kappa(x,\nu) < \max\left\{\dfrac{1}{k-1}, \dfrac{1}{2}\sqrt{\dfrac{n}{k-1}}\right\}$ *somewhere, or*

$$\kappa(x,\nu) \equiv \max\left\{\dfrac{1}{k-1}, \dfrac{1}{2}\sqrt{\dfrac{n}{k-1}}\right\}.$$

(iii) *For $M^2$ immersed in $S^4$, either $|\kappa(x,\nu)| < 1$ somewhere or $|\kappa(x,\nu)| \equiv 1$.*

**Proof.** Set $A = A(\mathcal{N}, \overline{\beta})$ to be an upper bound for $(2 \cdot \operatorname{vol}(\hat{N}(M)))^{2/n}$. The topological argument of [L-S, Cor. 4] shows that (5.7) must fail for $q = n$. Thus $\kappa \geq \max\left\{\dfrac{1}{k-1}, \dfrac{1}{2}\sqrt{\dfrac{n}{k-1}}\right\}$ everywhere implies $\kappa \equiv \max\left\{\dfrac{1}{k-1}, \dfrac{1}{2}\sqrt{\dfrac{n}{k-1}}\right\}$. This gives (ii), and (i) follows from [R-Y, Thm. 2.2]. A similar argument relying on [L-S, Cor. 6] gives (iii). □

These results improve [L-S], who show that for (ii), $\kappa \leq \max\left\{\dfrac{1}{k-1}, \dfrac{1}{2}\sqrt{\dfrac{n}{k-1}}\right\}$ somewhere, and for (iii), $|\kappa| \leq 1$ somewhere.

E. We now assume that $M$ is minimally isometrically immersed in $S^N$. The Gauss equation gives

$$1 = K(X,Y) - \langle \beta(X,X), \beta(Y,Y)\rangle + |\beta(X,Y)|^2 \tag{5.8}$$

for $X, Y \in T_x M$. Thus

$$\sum_{j=1}^{q} \sum_{\ell=q+1}^{n} \left(2|\beta(v_j, v_\ell)|^2 - \langle \beta(v_j, v_j), \beta(v_\ell, v_\ell)\rangle\right)$$

$$= \sum_{j=1}^{q} \sum_{\ell=q+1}^{n} \left(-K(v_j, v_\ell) + 1 + |\beta(v_j, v_\ell)|^2\right)$$

$$= q(n-q) + \sum_{j=1}^{q} \sum_{\ell=q+1}^{n} \left(|\beta(v_j, v_\ell)|^2 - K(v_j, v_\ell)\right)$$

$$= q(n-q) + \sum_{j=1}^{q} \sum_{\ell=q+1}^{n} \left(1 + \langle \beta(v_j, v_j), \beta(v_\ell, v_\ell)\rangle - 2K(v_j, v_\ell)\right)$$

$$= 2q(n-q) - 2\sum_{j=1}^{q} \sum_{\ell=q+1}^{n} K(v_j, v_\ell) + \langle \sum_{j=1}^{q} \beta(v_j, v_j), \sum_{\ell=q+1}^{n} \beta(v_\ell, v_\ell)\rangle$$

$$= 2q(n-q) - 2\sum_{j=1}^{q} \sum_{\ell=q+1}^{n} K(v_j, v_\ell) - |\sum_{j=1}^{q} \beta(v_j, v_j)|^2, \tag{5.9}$$

since $M$ minimal implies $\sum_{j=1}^{q} \beta(v_j, v_j) = - \sum_{\ell=q+1}^{n} \beta(v_\ell, v_\ell)$.

**Theorem 5E.** *Let $M^n$ be minimally isometrically immersed in $S^N$. If*

$$\Delta + \inf_{\substack{v_j, v_\ell \in T_x M \\ v_j \perp v_\ell}} \{\sum_{j=1}^{q} \sum_{\ell=q+1}^{n} K(v_j, v_\ell)\} - \frac{q(n-q)}{2} > 0,$$

*then $H_q(M; \mathbb{Z}) = H_{n-q}(M; \mathbb{Z}) = 0$. Thus if for all orthonormal bases $\{v_1, \ldots, v_n\}$ of $T_x M$,*

$$\sum_{j=1}^{q} \sum_{\ell=q+1}^{n} K(v_j, v_\ell) > \frac{q(n-q)}{2},$$

*then $H_q(M; \mathbb{Z}) = H_{n-q}(M; \mathbb{Z}) = 0$. In particular, if $K(v_j, v_\ell) > 1/2$ for all $v_j, v_\ell \in T_x M$, for all $x \in M$, then $M$ is a homotopy sphere. Finally, if $Ric > \frac{n-1}{2}$, then $\pi_1(M) = 0$, and if $\sum_{j=1}^{2} \sum_{\ell=3}^{n} K(v_j, v_\ell) > n - 2$, then $\pi_2(M) = 0$.*

The proof of the first statement follows from Theorem 5D and (5.9). For the last statement, we just note that $\sum_{\ell=2}^{n} K(v_1, v_\ell) = \langle \text{Ric}^\# v_1, v_1 \rangle$.

Observe that every such $M$ has $Ric \leq n - 1$; just set $X = v_1$, $Y = v_\ell$ in (5.8) and sum from $\ell = 2$ to $\ell = n$.

F. The same techniques give a generalization of the main results of [H-W].

**Theorem 5F.** *Let $M$ be isometrically immersed in $\mathbb{R}^N$ with*

$$\Delta + \inf_{\substack{v \in T_x M \\ |v|=1}} \{-\sum_i \|A(v, Z^i)\|^2 + \langle \text{Ric}^\# v, v \rangle\} > 0.$$

*Then*
*(i) for every compact Riemannian manifold $N$, there are no nonconstant stable harmonic maps $f: N \to M$, and the homotopy class of any map $\alpha_0: N \to M$ contains maps of arbitrarily small energy;*

*(ii) there are no nonconstant stable harmonic maps $f: M \to N$, and the homotopy class of any map $\alpha_0: M \to N$ contains maps of arbitrarily small energy;*

*(iii) $\pi_1(M) = \pi_2(M) = 0$.*

**Proof.** (i) and (ii) follow from Theorem 1B and Theorem 4B by setting $p = 2$, $q = 1$ in (4.2) and noting that $\delta^2 \Lambda^1 = 0$. Then (iii) follows from (i) as in Corollaries 1B.1, 1C.1. □

Howard and Wei prove (i)-(iii) under the assumption that the potential term satisfies
$$-\sum_i \|A(v, Z^i)\|^2 + \langle \text{Ric}^\# v, v \rangle > 0$$
for all non-zero $v \in T_x M$, for all $x \in M$. (To see that this condition agrees with [H-W], we just note that $\|A(v, Z^i)\|^2 = \langle v, A(A(v, Z^i), Z^i) \rangle = \langle v, A_{Z^i}^2 v \rangle$ in their notation, and use [H-W, (2.12)]. As in Theorem 5A, the hypothesis of Theorem 5F holds under a restriction on the $n/2$ norm of the potential.

G. We conclude with a few remarks and open questions.

1. The Laplacian $\Delta$ can be replaced by the Bismut-Witten Laplacian $\Delta^h = \Delta - 2\mathcal{L}_{\nabla h}$ for any $h \in C^\infty(M)$, where $\mathcal{L}$ denotes Lie bracket. As in [E-R2], [Li2], Theorem 5A generalizes to give vanishing of integral homology if

$$\Delta + R_0^q - \sup_{\substack{v \in \Lambda^q T_x M \\ |v|=1}} \{\sum_{j=1}^{q} \sum_{\ell=q+1}^{n} |\alpha(v_j, v_\ell)|^2\} + |\nabla h|^2 - \Delta H - 2(d\Lambda^q)(\text{Hess } h) > 0.$$

2. It would be nice to show that $\pi_2(G) = 0$ for any Lie group $G$ using the methods of this section. Since $G$ retracts onto its maximal compact subgroup, it follows from [Mi, §21] that we may assume $G$ is compact, simply connected, and simple. Let $G$ be given the biinvariant metric $g$ induced by the Killing form on the symmetric space $G \times G/\Delta G \simeq G$. It is well known that the (suitably normalized) eigenfunctions of the first eigenvalue $\lambda_1$ of the Laplacian give a minimal isometric immersion of $(G, \frac{\lambda_1}{n} g)$ into the unit sphere, where $n = \dim G$ (cf. [L]). By [O], for this immersion we have $\|\beta\|^2 = n(n-1) - (n^2/2\lambda_1)$. Using the standard inclusion $S^N \hookrightarrow \mathbb{R}^{N+1}$, Theorem 5A and Proposition 5B, we see that $\pi_2(G) = 0$ if

$$\frac{n}{\lambda_1} \sum_{j=1}^{2} \sum_{\ell=3}^{n} K(v_j, v_\ell) > \frac{n(n-1)}{2} - \frac{n^2}{4\lambda_1} + \frac{n}{2},$$

where $K$ denotes sectional curvature with respect to $g$.

Since $K(v_j, v_\ell) = \frac{1}{4}|[v_j, v_\ell]|^2$, $\pi_2(G) = 0$ if

$$\sum_{j=1}^{2} \sum_{\ell=3}^{n} |[v_j, v_\ell]|^2 > n(2\lambda_1 - 1) \tag{5.10}$$

for all orthonormal bases $\{v_i\}$ of the Lie algebra of $G$. The lowest eigenvalue $\lambda_1$ can be explicitly computed for the classical and exceptional simple groups [R]. However, it seems difficult to minimize the left hand side of (5.10) on a case by case basis. On the other hand, from the classical fact that $H^3(G;\mathbb{R}) \neq 0$ for any compact $G$, we see that $\lambda_1 > 1/2$ for such $G$. (Here we must note that the metric induced by the Killing form on $G$ as the symmetric space $G \times G/\Delta G$ is twice that induced on $G$ by its Killing form, so we have $\lambda_1 \geq 1/4$ for the metric on $G$ induced by the Killing form.)

By a similar computation, we see that $\pi_1(G) = 0$ if $\text{Ric}(v_1, v_1) > n(4\lambda_1 - 1)$ where $\lambda_1$ is computed for the Laplacian on $G$ (not on $G \times G/\Delta G$). Thus $\pi_1(G) \neq 0$ implies $\lambda_1 \geq \dfrac{1}{4} + \dfrac{\text{Ric}_0}{4n}$. This gives an improved estimate for $\lambda_1$ if $G$ has discrete center, as its sectional curvatures are then strictly positive.

3. Let $\pi\colon M' \to M$ be a $p$-fold cover of compact manifolds. Triangulate $M$ and lift the triangulation to $M'$. The map

$$\pi\colon \Sigma n_i \sigma_i \mapsto \Sigma n_i \pi^{-1}(\sigma_i)$$

taking simplicial chains on $M$ to simplicial chains on $M'$ induces a pullback map $\pi$ on homology.

**Lemma 5G.** $\pi\colon H_q(M;\mathbb{Z}\left[\frac{1}{p}\right]) \to H_q(M';\mathbb{Z}\left[\frac{1}{p}\right])$ *is injective.*

**Proof.** If $\alpha = \sum n_i \pi^{-1}(\sigma_i) = \partial \sum m_i \tilde{\sigma}_i$ for some simplices $\tilde{\sigma}_i$ in $M'$, then since $\alpha$ is invariant under deck transformations $\gamma \in \{\gamma^1, \ldots, \gamma^p\}$, $\alpha = \partial \sum_i m_i \gamma_* \tilde{\sigma}_i$. Thus

$$\alpha = \partial \sum_i \sum_{j=1}^p \frac{m_i}{p} \gamma_*^j \tilde{\sigma}_i = \partial \sum_i \frac{m_i}{p} \pi^{-1}(\sigma_i).$$

**Corollary 5G.** *Let $M$ satisfy the hypothesis of Theorem 5A. If $M \to N$ is a $p$-fold cover, then*

$$H_q(N;\mathbb{Z}\left[\frac{1}{p}\right]) = H_{n-q}(N;\mathbb{Z}\left[\frac{1}{p}\right]) = 0.$$

A typical example is $M = S^3$, $N = S^3/\mathbb{Z}_p$ a lens space, and $q = 1$. The case $q = 2$ is particularly intriguing, since $\pi_2(M)$ is unchanged for finite covers. If $\pi_1(M)$ is infinite and $R^2 > 0$, one might hope to find a finite cover $M'$ of $M$ with an isometric embedding of the pullback metric into $\mathbb{R}^N$ such that $\|\alpha\|^2$ for $M'$ is smaller than $\|\alpha\|^2$ for $M$. If by passing to a high enough cover $M'$, one can make $\|\alpha\|^2$ arbitrarily small,

then eventually $\pi_2(M') = 0$ and so $\pi_2(M) = 0$. (At present, it is known that $R^2 > 0$ implies $\pi_2(M)$ is a torsion group [E-R1].) However, we do not know how to construct such embeddings.

\*

# A   Appendix: Extrinsic formulae for Weitzenböck curvatures and the Hodge Laplacian

A.  The following is taken from [E3] but with the 'probabilistic' arguments rephrased. We consider an isometric immersion $f\colon M \to \mathbb{R}^m$ and use the notation of §4. In particular for $i = 1,\ldots,m$, $S_t^i\colon M \to M$, $t \in \mathbb{R}$, is the gradient flow for the $i^{\text{th}}$ component $f^i$ of $f$, and $X^i = \nabla f^i$.

**Proposition A.** *For any $C^2$ $q$-form $\psi$ on $M$*

$$\Delta^q \psi = -\sum_{i=1}^m \mathcal{L}_{X^i}\mathcal{L}_{X^i}\psi$$

*where $\mathcal{L}_{X^i}$ denotes Lie differentiation in the direction $X^i$. Moreover, for each $x \in M$ the Weitzenböck curvature $R_x^q$ has adjoint $(R_x^q)^*\colon \Lambda^q T_x M \to \Lambda^q T_x M$ given by*

$$(R_x^q)^*(V) = -\sum_{i=1}^m \frac{D^2}{\partial t^2}\Lambda^q(TS_t^i)(V)\bigg|_{t=0}.$$

**Proof.** Since $\sum_i \nabla X^i(X^i(x)) = 0$ for all $x \in M$, (because $\nabla X^i(x) = A_x(\cdot, e_i - X^i(x))$ for $A$ the shape operator and $e_1,\ldots,e_m$ the standard base for $\mathbb{R}^m$) the results are true for $q = 0$ (and well known). This vanishing together with a direct computation gives

$$\sum_{i=1}^m (\mathcal{L}_{X^i}\mathcal{L}_{X^i}\psi)_x = \text{trace } (\nabla^2 \psi)_x + \psi \circ Q_x^q \tag{A.1}$$

where $Q_x^q\colon \Lambda^q T_x M \to \Lambda^q T_x M$ is given by $Q_x^q(V) = \sum_i \frac{D^2}{\partial t^2}\Lambda^q(TS_t^i)(V)|_{t=0}$.

Suppose inductively that the claims are true for a fixed $q \in \{0, 1, \ldots, n-1\}$. Then, for a $C^2$ $q$-form $\psi$, since exterior differentiation and Lie differentiation commute, we

have

$$\begin{aligned}\Delta^{q+1}d\psi &= d\Delta^q\psi = -\sum_{i=1}^{m}\mathcal{L}_{X^i}\mathcal{L}_{X^i}d\psi \\ &= -\text{trace}\,(\nabla^2 d\psi) - (d\psi)\circ Q^{q+1}.\end{aligned}$$

But $\Delta^{q+1}d\psi = -\text{trace}\,(\nabla^2 d\psi)_x + (d\psi)\circ(R^{q+1})^*$ by the Weitzenböck formula. Thus $(d\psi)\circ Q_x^{q+1} = -(d\psi)\circ(R_x^{q+1})^*$ for all $x \in M$ and all $C^2$ $q$-forms $\psi$. Since any $\ell \in (\Lambda^{q+1}T_xM)^*$ has the form $\ell = (d\psi)_x$ for some $C^2$ $q$-form $\psi$, this shows $Q_x^{q+1} = -(R_x^{q+1})^*$ and from this follows the first claim by the Weitzenböck formula and (A.1).

B.  Note also the corollary that comes from (3.1) and was used to get (4.2), namely:

$$(R_x^q)^* = d\Lambda^q(\text{Ric}_x^{\#}) - \delta^2\Lambda^q(A(-, e_i - X^i(x))). \tag{A.2}$$